\documentclass{elsart}
\usepackage{cite}
\usepackage{verbatim}
\usepackage{amsmath}
\usepackage{amssymb}
\usepackage{amstext}
\usepackage{dcolumn}
\usepackage{epsfig}

\usepackage{amssymb}
\usepackage{amssymb, amsmath}
\usepackage{graphicx}
\usepackage{epstopdf}
\usepackage{subfig}
\usepackage[utf8]{inputenc}

\newtheorem{theorem}{Theorem}
\newtheorem{definition}{Definition}

\newtheorem{lemma}{Lemma}[section]

\date{}

\begin{document}
\begin{frontmatter}

\title{On the existence-uniqueness and computation of solution of a coupled PDE-ODE system with application to cardiac electric 
activity}
\author {Meena Pargaei} and
\ead{meenamenu15@gmail.com}
\author {B. V. Rathish Kumar\corauthref{cor}}
\corauth[cor]{corresponding author.}\ead{bvrk@iitk.ac.in}
\address {Department of Mathematics and Statistics,
 Indian Institute of Technology,\\ Kanpur 208016, India}

\begin{abstract}
 In this study, we consider a system of degenerate reaction-diffusion equations, which govern the 
 electric activity in the heart with a diffusion term modeling the potential in the surrounding tissue
and the nonlinear ionic model proposed by Morris $\&$ Lecar. The global existence of a solution is 
established based on regularization argument using Fedo-Galerkin/Compactness approach. The uniqueness
of the solution is shown based on Gronwell's Lemma upon some special treatment of nonlinear terms.
The system of the continuous space-time model is first reduced to a semi-discrete time-dependent 
system based on finite element formulation, and then the fully discrete system is derived using the 
Backward Euler time stepping scheme. The numerical solution obtained using FreeFem++ are presented.
\end{abstract}


\end{frontmatter}

\section{Introduction}
An electrocardiogram examines for the problems with the electrical activity of your heart. Bidomain
model\cite{Luca} is used for describing the cardiac electrical activity.
It consists of two PDEs coupled to a system of ODEs, describing the electrical activity of the heart.
The two PDEs describes the dynamics of intracellular and extracellular potentials,
whereas the ODEs, also known as an ionic model, describes the electrical behavior of the myocardium 
cell membrane. Involvement of different space and time scales makes it computationally expensive.

A simplified mathematical model of the cardiac tissue is the anisotropic monodomain system [7], 
which consist of a parabolic reaction-diffusion equation, describing the evolution of the membrane
potential, coupled with an ionic model. Some of the well known ionic models are FHN, Winfree, 
Rogers and McCulloch, Panfilov, Phase-I Luo Rudy, Morris Lecar model, etc.  

The first wellposedness of the Bidomain model with ionic model given by FHN
ionic model [5] has been proved in [2]. In [3] and [9], the existence of the solution is proved for 
a wide class of ionic models (including Panfilov [1] and McCulloch
[6]). Uniqueness, however, is achieved only for the FHN model. In [12] existence,
Uniqueness and some regularity results are proved with LR1 ionic model [4].

None of those mentioned above works consider the Morris Lecar ionic model. So, This paper describes
the existence of the Bidomain model with Morris Lecar as the ionic model. The main result states the
existence of a weak global solution for Bidomain equation with ionic model Morris Lecar. Numerical
simulation using the Finite Element Method for the monodomain model with Morris Lecar ionic model 
is also done here. For implemantaion FreeFem++ software is used.

In the next section, we describe the mathematical model (anisotropic Bidomain and monodomain models)
and ionic model also. In section 3 we state our main existence result for Bidomain equation with
Morris Lecar ionic model. In section 4 we give the proof of this result. Regularization argument 
and Fedo-Galerkin Technique is used to proof the result. In section 5 we provide some numerical 
results for the monodomain model with the Morris Lecar ionic model.

\section{Mathematical model}
\subsection{Bidomain model}
Bidomain model is a representation of the cardiac tissue as the superposition of two continuous 
anisotropic media, the intra (i) and extra (e) cellular media, coexisting at every point of the 
tissue and connected by a distributed continuous cellular membrane. This model describes the 
averaged intracellular and extracellular electric potentials and currents by a reaction-diffusion 
system of the degenerate parabolic type.

Let $\Omega\subset$ be the cardiac tissue domain.  For the  Bidomain characterization of the cardiac
tissue, it is considered as the overlapping of the intracellular and extracellular continuous domains
such that each point in the intracellular myocardium is also in the extracellular and the intracellular,
and the extracellular medium, is identified by the conductivity tensors $D_i$ and $D_e$. 
Let $a_l(x),a_t(x),a_n(x)$ be orthonormal  set of triplets  corresponding to the structure of the
cardiac tissue at a point x, where, $a_l$, parallel to the local fiber direction, and $a_n$, normal
to the cardiac muscle sheet.  Let ${\sigma_l}^{i,e}$,${\sigma_t}^{i,e}$, ${\sigma_n}^{i,e}$ are the 
conductivity coefficients along the corresponding directions. In general, these coefficients may 
depend on $x$, but in the following, we assume that they are constant, i.e., homogeneous anisotropy.
Then the conductivity tensors $D_i$ and $D_e$, generally dependent on the position $x$, is given by
\begin{equation}
\label{Die}
D_{i,e}(x) ={\sigma_l}^{i,e} a_l(x) {a_l}^T(x) + {\sigma_t}^{i,e} a_t(x) {a_t}^T(x) +{\sigma_n}^{i,e} a_n(x) {a_n}^T(x).
\end{equation}
If ${\sigma_n}^{i,e}$=${\sigma_t}^{i,e}$ we recover the axially isotropic case
\begin{align*}
D_{i,e}(x) = {\sigma_t}^{i,e} I +({\sigma_l}^{i,e}-{\sigma_t}^{i,e} )a_l(x) {a_l}^T(x).
\end{align*}
The bioelectric activity of cardiac cells is due to the flow $I_{ion}$ (per unit area of the membrane
surface) of various ionic currents (the most important being sodium, potassium, and calcium) through
the cellular membrane. Since the membrane behaves as a capacitor, the total membrane current per 
unit volume is given by  $I_m = \chi( C_m \frac{\partial v}{\partial t}+I_{ion})$ where $v=u_i-u_e$ 
is the transmembrane potential, the coefficient $\chi$ is the ratio of membrane area per tissue 
volume, $C_m$ is the surface capacitance of the membrane, and $I_{ion}$ is the ionic current 
described later and depending on the membrane model.

Imposing the conservation of currents, we have $divJ_i = - divJ_e = I_m$ where $J_{i,e}= -D_{i,e} u_{i,e}$, 
are the intracellualr and the extracellular current densities. Therefore, in the Bidomain model, 
the intra and extracellular potentials $u_i , u_e$ are modeled by the following reaction–diffusion 
system of PDEs, coupled with a system of ODEs for gating variables, descibed later. Given an applied
current per unit volume $I_{app} \colon\Omega \times (0,T) \rightarrow R$ , initial conditions
$v_0 \colon\Omega \rightarrow R$ , $w_0 \colon\Omega \rightarrow R^M$, find the intra and 
extracellular potentials $u_i,u_e \colon\Omega \rightarrow R$ , the transmembrane potential 
$v=u_i-u_e$ and the gating variables $w_0 \colon\Omega \rightarrow R^M$, such that
\begin{equation}
\label{ui}
C_m \frac{\partial v}{\partial t}- div(D_i(x)\nabla u_i) + I_{ion}(v,w)=I^i_{app} \hspace{1cm} in \Omega \times (0,T)
\end{equation}
\begin{equation}
\label{ue}
C_m \frac{\partial v}{\partial t}- div(D_e(x)\nabla u_e) + I_{ion}(v,w)=I^e_{app} \hspace{1cm} in \Omega \times (0,T)
\end{equation}
\begin{equation}
\label{w}
\frac{\partial w}{\partial t}-g(v,w)=0 \hspace{1cm} in  \Omega \times (0,T)
\end{equation}

We assume that the cardiac tissue is insulated, therefore homogeneous Neumann boundary conditions are assigned on $\partial \Omega \times (0,T)$ as $n^T D_{i,e} \nabla u_{i,e} =0$.

Initial conditions are assigned in $\Omega$ for $t = 0$ as follows
\begin{align*}
v(x,0)=u_i(x,0)-u_e(x,0) = v_0(x,0),~~ w(x,0)=w_0(x,0). 
\end{align*}
we then have the following compatibility condition for the system to be solvable:
\begin{align*}
\int_\Omega I^i_{app} = \int_\Omega I^e_{app}
\end{align*}

\subsection{Simplified monodomain model}

Assume that the anisotropy ratio of the two continuous media are equal, i.e. $D_i = \lambda D_e$ 
with $\lambda$ constant, and setting $D=\frac{\lambda D_i}{1+\lambda}$ and
$I_{app}=\frac{\lambda I^i_{app}}{1+\lambda} + \frac{I^e_{app}}{1+\lambda}$, then the Bidomain 
system reduces to the anisotropic Monodomain model consisting in a parabolic
reaction–diffusion equation for the transmembrane potential $v$ only which is descibed as by the following
set of equations:
\begin{equation}
\label{mv}
C_m \frac{\partial v}{\partial t}- div(D(x)\nabla v) + I_{ion}(v,w)=I_{app} \hspace{1cm} in \Omega \times (0,T)\\
\end{equation}
\begin{equation}
\label{mw}
\frac{\partial w}{\partial t}-g(v,w)=0 \hspace{1cm} in  \Omega \times (0,T)
\end{equation}

with Neumann boundary condition for $v$ and initial conditions for $v$ and $w$.

The conductivity tensor in the axial symmetric case is given by
\begin{align*}
D(x)= {\sigma_t} I +({\sigma_l}-{\sigma_t} )a_l(x) {a_l}^T(x)~~\text{with}~~\sigma_{l,t}=\frac{\lambda \sigma^i_{l,t}}{1+\lambda}.
\end{align*}
This model has been extremely used in computation because it requires substantially less computational and memory resources than the Bidomain
model. Nevertheless, it is not an adequate cardiac model since it is unable to reproduce some patterns and morphology of the experimentally observed extracellular
potential maps and electrograms. Therefore, unequal anisotropy ratio of the intra and extracellular media cannot be neglected.

\subsection{Ionic model}
\textbf {Morris Lecar Model\cite{ml}}

Morris Lecar model is a biological neuron model developed by Dr. Catherine Morris and Dr. Harold 
Lecar. A variety of oscillatory behavior of $Ca^{++}$ and $K^+$ conductance in the giant barnacle 
muscle fiber is replicated.

It consists of a two-dimensional system of non-linear ordinary differential equations. It is a 
simplified model version of the Hodgkin-Huxley ionic model. This system of equations qualitatively
describes the complex relationship between cellular membrane potential and the ion channel activation,
within the membrane: the potential depends on the activity of ion channels, the activity of ion 
channel depends on the voltage. Ionic current and the dynamics of gating variables are given by:
\begin{equation}
\label{Iion}
-I_{ion}(v,w) = \frac{1}{C_m} (g_{Ca}m_{\infty}(v)(v-v_{Ca})+ g_K w(v-v_K)+ g_L(v-v_L))
\end{equation}
\begin{equation}
\label{g}
g(v,w) = \phi \frac{w-w_{\infty}(v)}{\tau_w(v)}
\end{equation}
where

   $ m_{\infty}(v) = 0.5[1+tanh(\frac{v-v_1}{v_2})]$
   
    $w_{\infty}(v) = 0.5[1+tanh(\frac{v-v_3}{v_4})]$
    
    $\tau_w(v) = 1/cosh(\frac{v-v_3}{2v_4})$
    
    $v_L , v_{Ca} , v_K =$ equilibrium potential corresponding to leak , $Ca^{++} , K^+$ conductances, respectively.
    
    $v_1=$ potential at which $m_{\infty}=0.5$ 
    
    $v_2=$ reciprocal of slope of voltage dependence of $m_{\infty}$ 
    
    $v_3=$ potential at which $w_{\infty}=0.5$ 
    
    $v_4=$ reciprocal of slope of voltage dependence of $w_{\infty}$

\section{Existence of the bidomain model with Morris Lecar ionic model}
\textbf{Assumption}

\textbf{(A)} We Assume that the conductivities of the intracellular and the extracellular spaces
$D_i$, $D_e$ $\in [L_{\infty}(\Omega)]^{3\times3}$ are symmetric and uniformly positive definite,
i.e. there exist $\alpha_i>0 ,\alpha_e>0$ such that $\forall x \in R^3 , \forall \xi \in R^3$,\\
$\hspace{1cm}\xi^T D_i(x)\xi\geq \alpha_i {\mid\xi\mid}^2 \hspace{0.5cm} , \hspace{0.5cm}\xi^T D_e(x)\xi\geq \alpha_e {\mid\xi\mid}^2$

$\textbf{Notation:}$
$V=H^1(\Omega)= \{u:u\in L^2(\Omega), Du\in L^2(\Omega) \}$

\begin{definition}
\label{weak Formulation}
\textbf{(Weak Formulation)}
A weak solution of the Bidomain equation is a quadruplet of functions $(v,u_i,u_e,w)$ with the regularity 

$v\in L^{\infty}(0,T;H^1(\Omega))\cap H^1(0,T;L^2(\Omega))$

$u_i,u_e \in L^{\infty}(0,T;H^1(\Omega)) , w \in W^{1,\infty}(0,T;L^{\infty}(\Omega)) $ 

\begin{equation}
\label{weak ui}
C_m \int_{\Omega}\partial_tv \phi_i + \int_{\Omega}D_i \nabla u_i \nabla \phi_i + \int_{\Omega}I_{ion}(v,w)\phi_i = \int_{\Omega}I_{app}\phi_i
\end{equation}
\begin{equation}
\label{weak ue}
C_m \int_{\Omega}\partial_tv \phi_e - \int_{\Omega}D_e \nabla u_e \nabla \phi_e + \int_{\Omega}I_{ion}(v,w)\phi_e = \int_{\Omega}I_{app}\phi_e
\end{equation}
\begin{equation}
\label{weak w}
\partial_tw + g(v,w)=0
\end{equation}

for all $(\phi_i,\phi_e)\in H^1(\Omega) \times H^1(\Omega)$
Equations (\ref{weak ui}) and (\ref{weak ue}) holds in $D'(0,T)$ and equation (\ref{weak w}) holds a.e..
\end{definition}

The next theorem provides the main result for the existence of solution for the Bidomain equation
with Morris Lecar ionic model.

 \begin{theorem}
 \label{1}
  Let $T>0$, $I_{app} \in L^2(Q_T) , D_i , D_e \in [L^{\infty}(\Omega)]^{3\times 3}$ symmetric and 
  satisfying Assumption (A), $v_0 \in H^1(\Omega), w_0 \in L^{\infty}(\Omega)$ be the
 given data. If $w_0 \in L^{\infty}(\Omega)$ with a positive lower bound $r>0$, such that
 $r<w_0\leq 1$ in $\Omega$, then the problem ( \ref{ui}-\ref{w}) with (\ref{Iion}),
 (\ref{g}) and initial conditions)has a weak solution in the sense of Definition (\ref{weak Formulation}).
 \end{theorem}
The next section gives the proof of this theorem.

\section{Proof of the theorem 3.1}
The non-linear reaction-diffusion equations are degenerate in time. This issue is overcome here by 
adding a couple of Regularization terms, making bidomain equations parabolic. Regularization
and approximation of solution are merged here. Then the resulting regularized system can be analyzed
through a Fedo-Galerkin/compactness procedure and specific treatment of non-linear terms.

In section 4.1 , Regularization and Fedo-Galerkin Techniques are merged by introducing a regularized problem in finite dimension.

\subsection{A regularized problem in finite dimension}
Let ${\{h_k\}}_{k\in N^*}$ be a Hilbert basis of V. We assume that the basis functions are sufficiently
smooth and that ${\{h_k\}}_{k\in N^*}$ is an orthonormal basis in $L^2(\Omega)$.

For all $n\in N^*$, we define the finite-dimensional space $V_n$ generated by ${\{h_k\}}^n_{k=1}$ i.e. 
\begin{align*}
V_n=<{\{h_k\}}^n_{k=1}>
\end{align*}
Hence, we can introduce, for each $n\in N^*$, the following discrete problem associated with 
(\ref{weak ui}-\ref{weak w}): \textbf{Discrete Problem}
Find $(u_{i,n},u_{e,n})\in C^1(0,T;V_n\times V_n) , w_n\in C^1(0,T;L^{\infty}(\Omega))$ such that,
for $v_n=u_{i,n}-u_{e,n}$ and for all $(h,e)\in V_n\times V_n$ we have,
\begin{equation}
\label{disc ui}
C_m \int_{\Omega}\partial_tv_n h + \frac{1}{n}\int_{\Omega}\partial_tu_{i,n} h + \int_{\Omega}D_i \nabla u_{i,n} \nabla h + \int_{\Omega}I_{ion}(v_n,w_n)h
= \int_{\Omega}I_{app}h
\end{equation}
\begin{equation}
\label{disc ue}
C_m \int_{\Omega}\partial_tv_n e - \frac{1}{n}\int_{\Omega}\partial_tu_{e,n} e - \int_{\Omega}D_e \nabla u_{e,n} \nabla e + \int_{\Omega}I_{ion}(v_n,w_n)e
= \int_{\Omega}I_{app}e
\end{equation}
\begin{equation}
\label{disc w}
\partial_tw_n + g(v_n,w_n)=0  a.e. in Q_T 
\end{equation}

And verifying the initial conditions
\begin{equation}
v_n(0)=v_{0,n} , u_{i.n}(0)=u_{i,0,n} , u_{e.n}(0)=u_{e,0,n} ,
     w_n(0)=w_0   a.e. in \Omega 
\end{equation}
The auxiliary initial conditions for $u_{i,n}$ and  $u_{e,n}$ needed, are defined by introducing two arbitrary functions  $u_{i,0}$ ,  $u_{e,0}\in V$ such that
$v_0=u_{i,0}-u_{e,0}$ in $\Omega$. Then. for $n\in N^*$, we define $u_{i,0,n} , u_{e,0,n},w_{0,n}$ as the orthogonal projections on $V_n\times V_n\times V_n$, of
$u_{i,0} , u_{e,0},w_{0}$.

By construction of these sequences, we have
\begin{align*}
v_{0,n} , u_{i,0,n} , u_{e,0,n},w_{0,n} \rightarrow v_0 , u_{i,0} , u_{e,0},w_{0}~~\text{in}~~V^3\times L^2(\Omega).
\end{align*}

\subsection{Local existence of the discretized solution}
\begin{lemma}
\label{local existence}
Suppose that there exists $C$ such that

${\Vert u_{i,0,n}\Vert}_{H^1(\Omega)}$ + ${\Vert u_{e,0,n}\Vert}_{H^1(\Omega)}$ +${\Vert w_{0,n}\Vert}_{L^2(\Omega)}\leq C$

For all $n\in N^*$ there exists a positive time $0<t_n<T$ which only depends on C such that the Discrete problem (12)-(15)
admits a unique solution over the time interval $[0,t_n]$.
\end{lemma}
\textbf{Proof.} Since ${\{h_k\}}^n_{k=1}$ is a Hilbert Basis of $V_n$ and orthonormal basis of $L^2(\Omega)$.
Since $L^{\infty}(\Omega)\subset L^2(\Omega)$ so $w_0,w\in L^2(\Omega)$.
Now we are choosing $\theta \in V_n$ and multiplying on both side of equation (14), we get
\begin{equation}
\int_{\Omega}\partial_tw_n + \int_{\Omega}g(v_n,w_n)=0
\end{equation}
Since ${\{h_k\}}^n_{k=1}$ is a Hilbert Basis of $V_n$ and orthonormal basis of $L^2(\Omega)$.
So, now we can write,

$u_{i,n}=\sum^n_{l=1}c_{i,l}(t)h_l $ , $u_{e,n}=\sum^n_{l=1}c_{e,l}(t)h_l $ , $w_{n}=\sum^n_{l=1}c_{w,l}(t)h_l $ ,

$u_{i,0,n}=\sum^n_{l=1}c^0_{i,l}h_l $ , $u_{e,0,n}=\sum^n_{l=1}c^0_{e,l}h_l $ , $w_{0,n}=\sum^n_{l=1}c^0_{w,l}h_l $

\textbf{Notation}

$c_i=\{ c_{i,l}\}^n_{l=1}$ , $c_e=\{ c_{e,l}\}^n_{l=1}$ , $c_w=\{ c_{w,l}\}^n_{l=1}$

$c^0_i=\{ c^0_{i,l}\}^n_{l=1}$ , $c^0_e=\{ c^0_{e,l}\}^n_{l=1}$ , $c^0_w=\{ c^0_{w,l}\}^n_{l=1}$

Then the system of equations (12),(13),(16) is equivalent to the following non-linear system of ODEs

\begin{equation}
\label{ode system}
M\begin{bmatrix} c^{'}_i \\  c^{'}_e \\ c^{'}_w \end{bmatrix} = \begin{bmatrix} G_i(t,c_i,c_e,c_w) \\  G_e(t,c_i,c_e,c_w)\\ G_w(t,c_i,c_e,c_w)\end{bmatrix}
\end{equation}

$\begin{bmatrix} c_i(0) \\  c_e(0) \\ c_w(0) \end{bmatrix}$ = $\begin{bmatrix} c^{0}_i \\  c^{0}_e \\ c^{0}_w \end{bmatrix}$

Here the matrix $M\in R^{3n\times 3n}$ is given by

$M$= $\begin{bmatrix} (C_m+\frac{1}{n})M_V & \vdots-C_m M_V & \vdots & 0 \\ \cdots \hspace{5mm} \vdots & \cdots &\cdots \\ -C_m M_V
&(C_m+\frac{1}{n}) M_V & \vdots & 0\\ \cdots & \cdots & \cdots \\ 0 \hspace{1cm}\vdots & 0 \hspace{1cm}\vdots & M_V \end{bmatrix}$
\vspace{3mm}\\
with $M_V\in R^{n\times n}$

and $M_V= (\int_{\Omega}h_k h_l)_{1\leq k,l\leq n}$

and right-hand side of (\ref{ode system}) is given by

$G_i(t,c_i,c_e,c_w)= - \int_{\Omega}D_i \nabla u_{i,n} \nabla h_k - \int_{\Omega}I_{ion}(v_n,w_n)h_k + \int_{\Omega}I_{app}h_k$

for all $1\leq k \leq n$,

$G_i(t,c_i,c_e,c_w)= - \int_{\Omega}D_i \nabla u_{i,n} \nabla h_k - \int_{\Omega}I_{ion}(v_n,w_n)h_k + \int_{\Omega}I_{app}h_k$

for all $1\leq k \leq n$,

$G_w(t,c_i,c_e,c_w)= - \int_{\Omega}g(v_n,w_n)h_k$

for all $1\leq k \leq n$,

\begin{lemma}
\label{matrix pd}
 For all $n \in N^*$, the matrix M is positive definite.
\end{lemma}
According to Lemma (\ref{matrix pd}) the mass matrix M is positive definite and hence invertible and, on the other hand, the RHS of (\ref{ode system}) is a $C^1$ function with respect
to the arguments $c_i , c_e , c_w$. So, By using the Cauchy-Lipschitz theorem, existence of the local solution of the ODE system (\ref{ode system}) follows.

\subsection{Energy estimates}
\begin{lemma}
\label{energy estimate}
 Let $u_{i,0} , u_{e,0}\in H^1(\Omega)$ and $w_0 \in L^{\infty}(\Omega)$ with $r<w_0\leq 1$ , there exists a positive constant $w_{min}$ independent of $T'$
such that a solution $(u_{i,n},u_{e,n},w_n)$ of the discrete problem defined on $[0,T']$ for $T'>0$ satisfies

$\Vert v_n\Vert_{L^{\infty}(0,t;L^2(\Omega))} + \frac{1}{\sqrt n}(\Vert u_{i,n}\Vert_{L^{\infty}(0,t;L^2(\Omega))}
+ \Vert u_{e,n}\Vert_{L^{\infty}(0,t;L^2(\Omega))})
+\Vert\nabla u_{i,n} \Vert_{L^2(Q_t)} + \Vert \nabla u_{e,n} \Vert_{L^2(Q_t)} \leq c$

$\Vert \partial_t v_{n} \Vert_{L^2(Q_t)} + \Vert v_n\Vert_{L^{\infty}(0,t;H^1(\Omega))}
+ \frac{1}{\sqrt n}(\Vert \partial_t u_{i,n} \Vert_{L^2(Q_t)} + \Vert \partial_t u_{e,n} \Vert_{L^2(Q_t)})
+ \Vert \nabla u_{i,n}\Vert_{L^{\infty}(0,t;L^2(\Omega))} + \Vert \nabla u_{e,n}\Vert_{L^{\infty}(0,t;L^2(\Omega))}\leq c$

and, for all $t\in [0,T']$

$\Vert w_n\Vert_{W^{1, \infty}(0,t;L^{\infty}(\Omega))}\leq c$ , $w_{min}\leq w_n\leq 1$ in $Q_T$.
\end{lemma}
\textbf{Proof.}From equation (\ref{disc w}) it follows that

$\partial_tw_n = -g(v_n,w_n)$

$\partial_tw_n = -\phi \frac{w_n-w_{\infty}(v_n)}{\tau(v_n)}$ where  $0\leq w_{\infty}\leq 1$

$\partial_tw_n \geq -\phi \frac{w_n}{\tau(v_n)}$

$ w_n\geq w_0 exp(-\frac{T\phi}{\tau(v_n)})$

again

$\partial_tw_n \leq -\phi \frac{w_n-1}{\tau(v_n)} = \phi \frac{1-w_n}{\tau(v_n)} $

$ w_n\leq 1-(1-w_0) exp(-\frac{T\phi}{\tau(v_n)})$

using $r<w_0\leq 1$, we obtain that

$w_{min}= r exp(-\frac{T\phi}{\tau(v_n)})\leq w_n \leq 1$, a.e. in $Q_T$.

On the other hand,

$-\frac{\phi}{\tau_w}\leq \partial_t w_n\leq \frac{\phi}{\tau_w} $

$\Vert w_n\Vert_{W^{1, \infty}(0,t;L^{\infty}(\Omega))}\leq \frac{\phi}{\tau_w}$

For first energy estimate take $h=u_{i,n}$ and $e=u_{e,n}$ in equation (\ref{disc ui}) and (\ref{disc ue}) and then subtract, which yields

$\frac{1}{2}\frac{d}{dt}[C_m \Vert v_n \Vert^2_{L^2(\Omega)} + \frac{1}{n} (\Vert u_{i,n} \Vert^2_{L^2(\Omega)}+ \Vert u_{e,n} \Vert^2_{L^2(\Omega)})]
+\alpha_i \Vert \nabla u_{i,n} \Vert^2_{L^2(\Omega)}$

$+ \alpha_e \Vert \nabla u_{e,n} \Vert^2_{L^2(\Omega)}
+ \int_{\Omega}I_{ion}(v_n,w_n)v_n \leq \int_{\Omega}I_{app}v_n$

here,

$-I_{ion}(v,w) = \frac{1}{C_m} (g_{Ca}m_{\infty}(v)(v-v_{Ca})+ g_K w(v-v_K)+ g_L(v-v_L))$

$I_{ion}v = -\alpha v^2-\beta v^2 w-\gamma v-\delta vw$

where, $\alpha = \frac{1}{C_m} (g_{Ca}m_{\infty}(v)+g_L)$ , $\beta=\frac{1}{C_m}g_K$ ,

$\gamma = \frac{1}{C_m}(v_{Ca}g_{Ca}m_{\infty}(v)+g_Lv_L)$ , $\delta = -\frac{1}{C_m}g_K v_K$

So, $I_{ion}v\geq -\alpha^{'}\mid v \mid^2-\gamma^{'}$ (because $\alpha$ and $\gamma$ are bounded)

where $\alpha^{'}$ and $\gamma^{'}$ are positive constants.

Therefore,

$\frac{1}{2}\frac{d}{dt}[C_m \Vert v_n \Vert^2_{L^2(\Omega)} + \frac{1}{n} (\Vert u_{i,n} \Vert^2_{L^2(\Omega)}
+ \Vert u_{e,n} \Vert^2_{L^2(\Omega)})]
+\alpha_i \Vert \nabla u_{i,n} \Vert^2_{L^2(\Omega)}$

$+ \alpha_e \Vert \nabla u_{e,n} \Vert^2_{L^2(\Omega)} \leq (\alpha^{'}+\frac{1}{2})\Vert v_n \Vert^2_{L^2(\Omega)}
+ \frac{1}{2})\Vert I_{app} \Vert^2_{L^2(\Omega)} +\gamma^{'}\mid \Omega \mid $

Therefore, integrating over $(0,t)$ with $t\in [0,T']$, and after that applying Gronwall Lemma and using the fact that,

$\frac{1}{n}({\Vert u_{i,0,n}\Vert}_{L^2(\Omega)} + {\Vert u_{e,0,n}\Vert}_{L^2(\Omega)})$ +${\Vert v_{0,n}\Vert}^2_{L^2(\Omega)}$ is uniformly bounded with
respect to $n$, we will get first estimate.

For the estimate of the time derivatives, we take $h= \partial_t u_{i,n}$ and $e=\partial_t u_{e,n}$ in equation (\ref{disc ui}) and (\ref{disc ue})
and subtract then integrating over $(0,t)$
with $t\in [0,T']$, we obtain

$\frac{1}{4}C_m\Vert \partial_t v_{n} \Vert^2_{L^2(Q_t)}
+ \frac{1}{n}(\Vert \partial_t u_{i,n} \Vert^2_{L^2(Q_t)} + \Vert \partial_t u_{e,n} \Vert^2_{L^2(Q_t)})
+ \frac{\alpha_i}{2}\Vert \nabla u_{i,n}\Vert_{L^2(\Omega)} + \Vert \nabla u_{e,n}\Vert_{L^2(\Omega)}$

$\leq c(\Vert \nabla u_{i,0,n} \Vert_{L^2()\Omega})+ \Vert \nabla u_{e,0,n} \Vert_{L^2()\Omega}) + \frac{1}{2} \Vert I_{app}\Vert_{L^2(Q_t)}
+ \gamma^{'}T+\alpha^{'}\Vert v_n \Vert^2_{L^2(Q_t)} + \int^t_0\int_{\Omega}(\beta v_n+\delta)w_n\partial_tv_n$

The last term is solved as
$\int^t_0\int_{\Omega}(\beta v_n+\delta)w_n\partial_tv_n \leq \int^t_0\int_{\Omega}(\beta v_n+\delta)\partial_tv_n \leq \Vert
(\beta v_n+\delta)\Vert_{L^2(Q_t)} \Vert \partial_t v_n \Vert_{L^2(Q_t)}$.

So, by inserting this inequality in the last expression and using the previous estimate we obtain the second estimate.

\subsection{Global existence of solution}
Energy estimates allows us to extend the existence time of our discrete solution $u_{i,n},u_{e,n},w_n$.
According to Lemma (\ref{energy estimate}),the solution satisfies, for all $t\in [0,T']$, where $T'$ is the existence time

${\Vert u_{i,n}(t)\Vert}_{H^1(\Omega)}$ + ${\Vert u_{e,n}(t)\Vert}_{H^1(\Omega)}$ +${\Vert w_{n}(t)\Vert}_{L^2(\Omega)}\leq c$

After applying Lemma (\ref{local existence}) iteratively, we obtain the existence of solution upto an arbitrary time $T$.

We now want to pass to the limit when n goes to infinity. Let us multiply equation (\ref{disc ui}) (\ref{disc ue}) by a function $\eta \in D(0,T)$ and integrate between $0 \& T$.
For all $k\leq n$, we have

$C_m \int^T_0 \int_{\Omega}\eta \partial_tv_n h_k + \frac{1}{n} \int^T_0  \int_{\Omega}\eta \partial_tu_{i,n} h_k + \int^T_0  \int_{\Omega}\eta D_i \nabla u_{i,n} \nabla h_k $\\
$+ \int^T_0  \int_{\Omega}\eta I_{ion}(v_n,w_n)h_k
= \int^T_0  \int_{\Omega}\eta I_{app}h$

$C_m \int^T_0 \int_{\Omega}\eta \partial_t v_n h_k - \frac{1}{n} \int^T_0  \int_{\Omega}\eta \partial_tu_{e,n} h_k - \int^T_0  \int_{\Omega}\eta D_e \nabla u_{e,n} \nabla h_k $

$+ \int^T_0  \int_{\Omega}\eta I_{ion}(v_n,w_n)h_k
= \int^T_0  \int_{\Omega}\eta I_{app}h$

From Lemma (\ref{energy estimate}) it follows that there exists $u_i , u_e \in L^{\infty}(0,T;H^1(\Omega)) , v \in L^{\infty}(0,T;H^1(\Omega))\cap H^1(0,T;L^2(\Omega)) , w \in L^{\infty}(Q_T)$ such that

$u_{i,e,n} \rightarrow u_{i,e}$ in $L^{\infty}(0,T;H^1(\Omega))$ weak $*$

$v_{n} \rightarrow v$ in $L^{\infty}(0,T;H^1(\Omega))$ weak $*$

$v_{n} \rightarrow v $ in $H^1(0,T;L^2(\Omega))$ weak

$w_n \rightarrow$ w in $L^{\infty}(Q_T)$ weak $*$

According to Lemma (\ref{energy estimate}), we also conclude that $\frac{1}{\sqrt n}u_{i,n}$ and $\frac{1}{\sqrt n}u_{e,n}$ are bounded in $L^{\infty}(0,T;L^2(\Omega))$.
So, for all $k \in N^*$
and $\alpha \in D(0,T)$, we have

$ lim_{n\to\infty} \frac{1}{n} \int^T_0  \int_{\Omega}\eta \partial_tu_{i,n} h_k = 0$ ,
$ lim_{n\to\infty} \frac{1}{n} \int^T_0  \int_{\Omega}\eta \partial_tu_{e,n} h_k = 0$

Since $v_n$ is bounded in $L^2(0,T;H^1(\Omega))\cap H^1(0,T;L^2(\Omega))$, so $v_n$ will be bounded in $H^1(Q_T)$. Hence by compact embedding of $H^1(Q_T)$ in $L^3(Q_T)$, the sequence
$\{v_n\}$ strongly converges to $v$ in $L^3(Q_T)$,

$ lim_{n\to\infty}\int^T_0  \int_{\Omega}\eta (\alpha(v_n) v_n+\gamma(v_n))h_k= \int^T_0  \int_{\Omega}\eta(\alpha(v) v+ \gamma(v)) h_k$

Since sequence $\{w_n\}$ is bounded in $L^{\infty}(Q_T)$ and $v_n$ strongly converges to $v$ in $L^2(Q_T)$, we have

$ lim_{n\to\infty}\int^T_0  \int_{\Omega}\eta (\beta v_n+\delta)w_n h_k= \int^T_0  \int_{\Omega}\eta (\beta v+\delta)w h_k$

Thus,
$ lim_{n\to\infty}\int^T_0  \int_{\Omega}\eta I_{iom}(v_n,w_n)h_k = \int^T_0  \int_{\Omega}\eta I_{ion}(v,w)h_k$

Since

$w_n \rightarrow$ w in $L^{\infty}(Q_T)$ weak $*$

\subsection{Uniqueness of the weak solution}
\begin{lemma}
\label{uniqueness}
Assume that the first partial derivatives of $I_{ion}(v,w)$ and $g(v,w)$ are bounded and that $(v_1,u_{i,1},u_{e,1},w_1) , (v_2,u_{i,2},u_{e,2},w_2)$ are two weak
solutions of our problem corresponding, respectively , to the initial data $(v_{1,0},w_{1,0})$ and $(v_{2,0},w_{2,0})$ and right-hand sides $I_{app,1}$ and $I_{app,1}$.
For all $t \in (0,T)$, there holds

$\Vert v_1(t)-v_2(t) \Vert^2_{L^2(\Omega)} + \Vert w_1(t)-w_2(t) \Vert^2_{L^2(\Omega)}$

$\leq exp(K_1t)K_2(\Vert v_{1,0}-v_{2,0} \Vert^2_{L^2(\Omega)} + \Vert w_{1,0}-w_{2,0} \Vert^2_{L^2(\Omega)}+ \Vert I_{app,1}-I_{app,2} \Vert^2_{L^2(Q_t)})$
\end{lemma}
\textbf{Proof:} Ref\cite{Boulakia}.
\indent

\textbf{Remark:} This result also provides a stability estimate with respect to the initial condition.

\section{Numerical method}
We consider Monodomain model with Morris Lecar ionic model and consider $I_{app}=0$. Our monodomain system is equivalent to finding $v\in H^1(\Omega)$ and $w\in L^{\infty}(\Omega)$ such that
\begin{equation}
C_m \int_{\Omega}\frac{\partial v}{\partial t}\zeta- \int_{\Omega}D(x)\nabla v \nabla \zeta + \int_{\Omega}I_{ion}(v,w)\zeta = 0 \\
\end{equation}
\begin{equation}
\int_{\Omega}\frac{\partial w}{\partial t}\varsigma -\int_{\Omega} g(v,w) \varsigma=0
\end{equation}

for all $\zeta \in H^1(\Omega) , \varsigma \in L^2(\Omega)$.
The system is discretized in space using finite element method and in time using Backward Euler method.

\subsection{Finite element discretization in space}
We consider the square domain $\Omega=[a,b]^2$. The domain$\Omega$ is discretized by introducing
a structured quasi-uniform grid of triangular $P1$ elements(which is denoted by $\tau_h$). So, FEM approximation for the domain $\Omega$ is, $\Omega = \cup_{E\in\tau_h}E$.

The associated finite element space 

$V_h = \{ \zeta_h \in V : \zeta_h \text{is continuous in} \Omega :
{\zeta_h}_{\mid E} \in P_1({E}), \forall E \in \tau_h \}$ 

A semi-discrete form is obtained by applying a standard Galerkin procedure and choosing a finite element basis $\{\zeta_i\} \in V_h$. Let $I^h_{ion}$ be the finite element approximation of
$I_{ion}$.

In the monodomain model, the finite element approximation $v_h$ of transmembrane potential $v$ and $w_h$ of gating variable $w$ are the solution of

\begin{equation}
 M \frac{\partial v_h}{\partial t}+ A v_h + MI^h_{ion}(v_h,w_h)=0 ,
 \frac{\partial w_h}{\partial t}=g(v_h,w_h)
\end{equation}
where

$M=(m_{rs}) , m_{rs}=\sum_E \int_{E}\zeta_r \zeta_s dx$

$A=(a_{rs}) , a_{rs}=\sum_E \int_{E}{\nabla \zeta_r}^T D(x) \nabla \zeta_s dx$

Numerical quadrature in 2-dimension is used in order to compute these integrals.
Now this ODE system is discretized by Backward Euler and the implementation is done using FreeFem++ library functions\cite{freefem}.

\subsection{Numerical results and discussion}

\begin{figure}[h]
\includegraphics[width=1\textwidth]{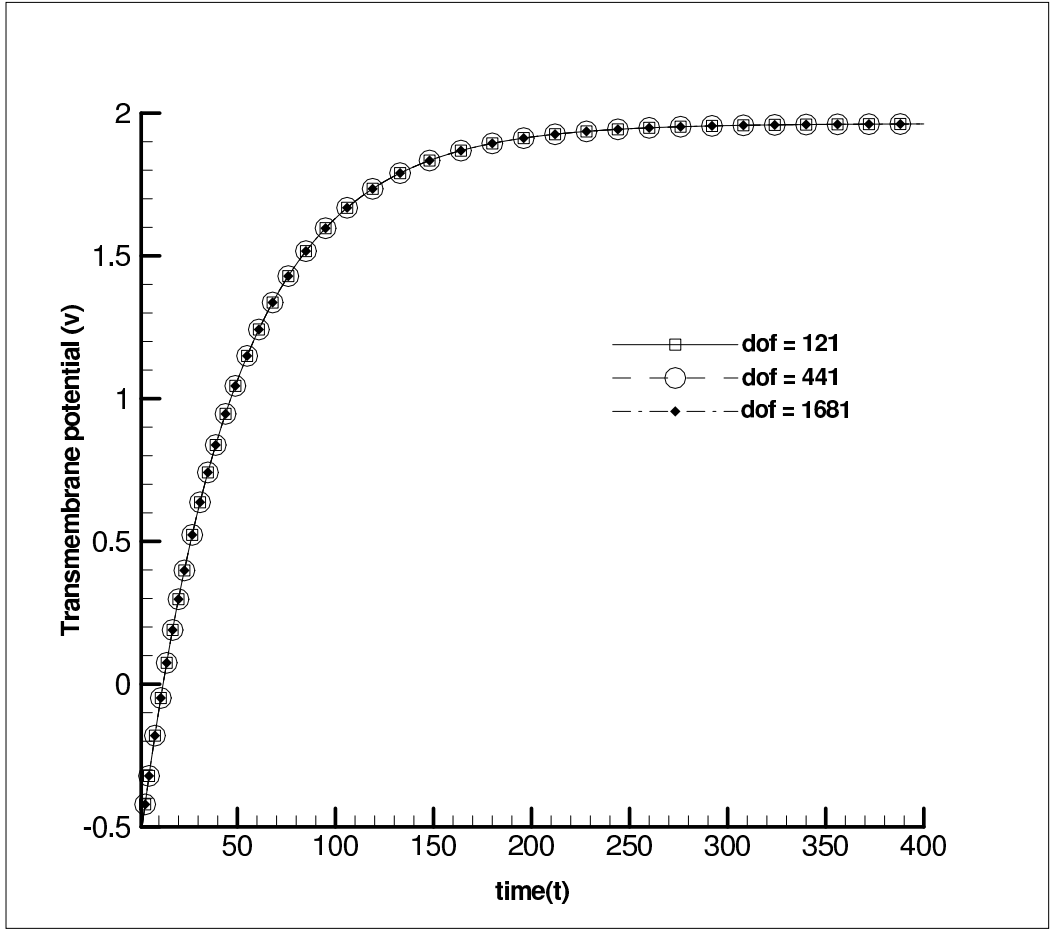}
\caption{plot corresponding to point (0,0) for dofs=121 , 441 , 1681.}
\label{t00-plot}
\end{figure}

\begin{figure}[h]
\includegraphics[width=1\textwidth]{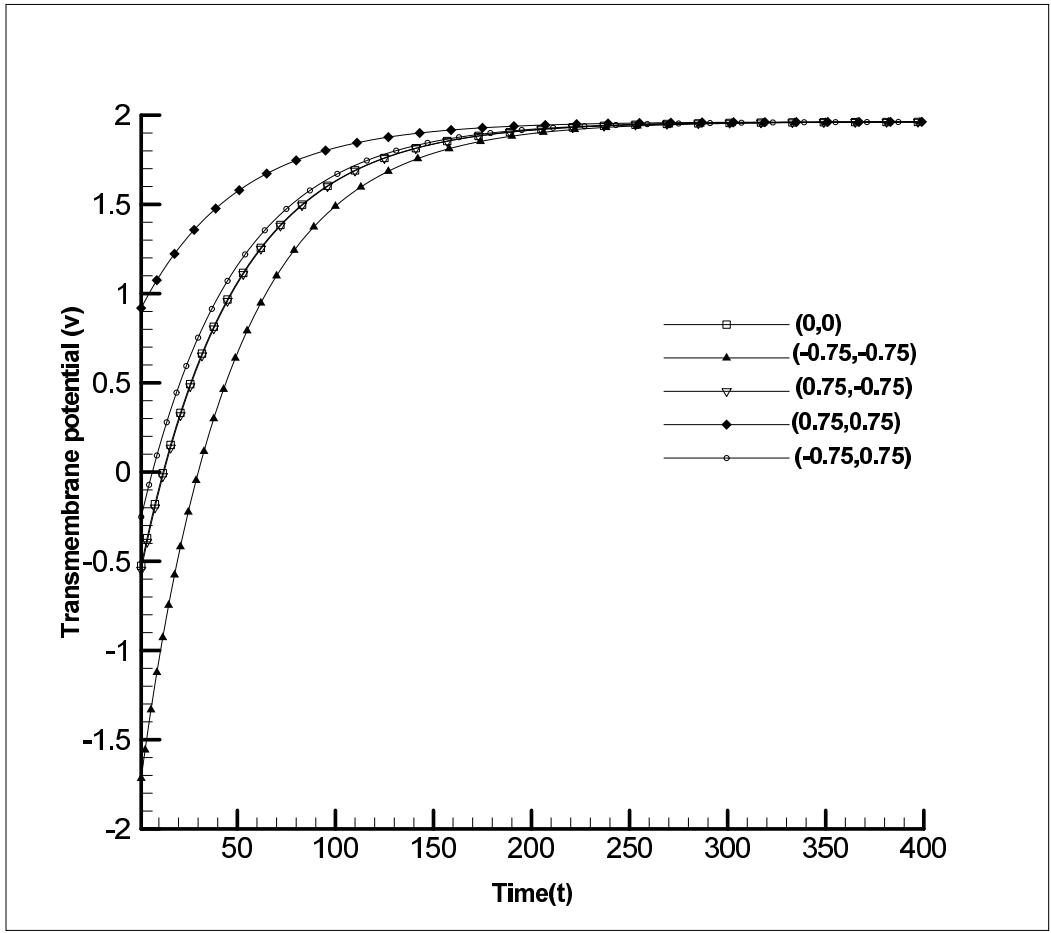}
\caption{plot for five different points with dof=1681.}
\label{test2}
\end{figure}

\begin{figure}[h]
\includegraphics[width= 1\textwidth]{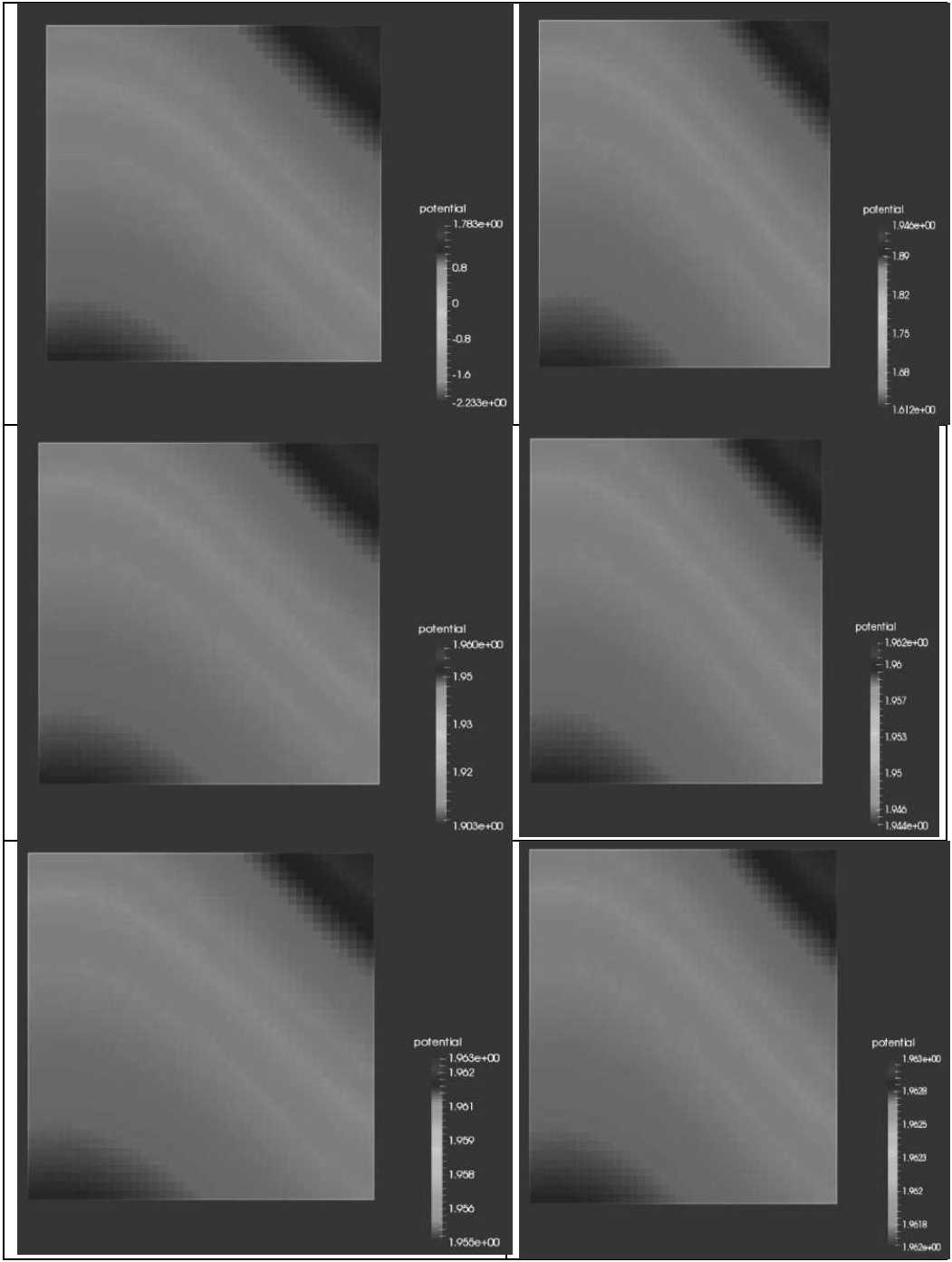}
\caption{the isochrones or iso-potential plots for the entire domain corresponding to six different times t=0 , 120 , 210 , 270 , 310 , 399 respectively.}
 \label{Fig:contour}
\end{figure}

To numerically simulate the transmembrane potential ($v$) in a cardiac tissue as modeled by the 
coupled PDE-ODE system (\ref{mv}-\ref{mw}) representing the monodomain model\cite{Luca} for cardiac
tissue together with the ionic model proposed by Morris and Lecar\cite{ml} the required parameter 
is chosen as given in\cite{Luca} and\cite{ml} and these details are provided in the table
\ref{tab:table1} and table \ref{tab:table2}.
In all the numerical simulations cardiac tissue has been represented by a square domain 
$\Omega=[-1.25,1.25]^2 $. To carry out the numerical simulations one has to choose an appropriate
grid system so that numerical solution is computed to the acceptable accuracy. Such a goal has been 
achieved here through grid validation test. Three different grid systems consisting of (a) 121, (b) 441 and (c) 1681 degrees of
freedom (dofs) ( or (a) 200, (b) 800 and (c) 3200 number of linear triangular elements respectively) 
have considered. Transmembrane potential ($v$) obtained using these three grid systems
are compared at different points of the domain.  In all these cases only a marginal variation 
(less than $0.5 \%$) in “$v$” is noticed as one moves from the grid system (a) to the grid system (c).
As a sample in Fig. \ref{t00-plot} the transmembrane potential corresponding to (x,y)=(0,0) with 
time are compared. Clearly, the grid system with 1681 dofs is more than adequate for the current set
of simulations. Hence all the simulations are carried out using 3200 linear triangular elements with
1681 dofs. In Fig. \ref{test2} the temporal variation of transmembrane potential corresponding
to the following five different representative points, chosen from the four different quadrants, of
the domain are presented. The initial transmembrane potential at these points follow the IC
setup, and they gradually evolve to the same steady state with $v \cong 2$, indicating that the 
intracellular and extracellular potentials have reached a state which no more favor calcium and 
Potassium ion migration and thereby maintain stable ionic concentration in the absence of any
depolarization phase. It is noticed that in about 300-time steps the transmembrane potential
corresponding to all these points already reach the steady-state and on further time marching only
a marginal variation in v is noticed. So the entire domain is nearly re-polarized in about 300-time
steps where each time step corresponds to 0.1-millisecond size. It is also to be noted lesser the 
initial transmembrane potential than its steady state value, more rapid is the growth in $ v$,
due to larger driving force (potential difference), and hence all points irrespective of its 
starting value would reach a steady state in about 300 time steps. This also indicates that
smaller the $v$ than 2, stronger may be the migration of calcium and potassium ions.  Now to 
trace the transmembrane potential in the entire domain in the form of isochrones or iso-potential
plots for the entire domain corresponding to six different time instances covering the initial state
to the steady state situation are presented in Fig. \ref{Fig:contour}. From the legend values for
$v$ it is clear that the min-max differences in $v$ reach to a nearly zero state in about 400 
time-steps each of 0.1 milliseconds.

\section{Conclusions}
We prove the Existence uniqueness for the solution of the bidomain model with Morris Lecar ionic 
model using Galerkin and Compactness approach. Also, Finite Element Computations based on Monodomain
model depict the success and effectiveness of Morris Lecar ionic model in enabling the visualization
the Cardiac Electric Activity in cardiac tissue.

\section{Acknowledgement}
We would like to thank the DST for the support through Inspire Fellowship. Also thankful to
FreeFem++ open source developers\cite{website} for making it available for the researchers.

\newpage
\begin{table}[ht]
  \centering
\caption{parameter Values used for the simulation}
\label{tab:table1}
  \begin{tabular}{c c c c c c c }
   parameter & $v_1$ & $v_2$ & $v_3$ & $v_4$ & $v_{Ca}$ & $v_K$  \\
    \hline
   values & -1.2 & 18 & -1 & 14.5 & 120 & -70   \\
   \end{tabular}
\end{table}

\begin{table}[ht]
  \centering
  \caption{parameter Values used for the simulation}

\label{tab:table2}
  \begin{tabular}{c c c c c c c c}
   parameter & $v_L$ & $g_l$ & $g_K$ & $g_{Ca}$ & $\lambda$ & $\sigma_l$ & $\sigma_t$\\
   \hline
   values  & -50 & 4 & 8 & 3 & 1 & $1.2*10^{(-3)}$ & $2.5562*10^{(-4)}$\\

\end{tabular}
\end{table}

%


\end{document}